\documentclass[amsmath,twocolumn,secnumarabic,superscriptaddress,%
    floatfix,amssymb,aps]{revtex4}
\usepackage{times,graphicx,amsthm,url}
\usepackage[percent]{overpic}
\newcommand{\nicefrac}[2]{\mbox{\kern.05em\raise.5ex\hbox{\footnotesize$#1$}
  \kern-.35em$/$\kern-.12em\lower.5ex\hbox{\footnotesize$#2$}\kern.04em}}

\graphicspath{{./figs/}}
\newcommand{\newbd}{\nicefrac{5\pi}{3}}
\newcommand{\upbd}{7.16}

\newcommand{\R}{\mathbb{R}}

\newcommand{\setm}{\smallsetminus}

\newcommand{\arc}[1]{\gamma_{#1}}
\newcommand{\len}[1]{\ell_{#1}}  
\newcommand{\lenc}[1]{\ell(#1)}  

\newcommand{\eps}{\varepsilon}
\renewcommand{\phi}{\varphi}

\newcommand{\0}{\mathbf{0}}

\newtheorem{theorem}{Theorem}[section]
\newtheorem{lemma}[theorem]{Lemma}
\newtheorem{corollary}[theorem]{Corollary}
\newtheorem{proposition}[theorem]{Proposition}

\theoremstyle{definition}
\newtheorem*{remark}{Remark}
\newtheorem*{definition}{Definition}

\newcommand{\figr}[1]{Figure~\ref{fig:#1}}
\newcommand{\thm}[1]{Theorem~\ref{thm:#1}}
\newcommand{\prop}[1]{Proposition~\ref{prop:#1}}
\newcommand{\cor}[1]{Corollary~\ref{cor:#1}}
\newcommand{\lem}[1]{Lemma~\ref{lem:#1}}

\newcommand{\arcsec}{\mathop{\textrm{arcsec}}}

\let\mgp=\marginpar \marginparwidth16mm \marginparsep1mm
\def\marginpar#1{\mgp{\raggedright\tiny #1}}
\def\marginpar#1{}   

\let\lbl=\label
\def\label#1{\lbl{#1}\ifinner\else\marginpar{\ref{#1}\ \\#1}\ignorespaces\fi}

\makeatletter
\newbox\overbox
\def\fakeover#1{\setbox\overbox\hbox{$#1$}\hbox
                         {$\overline{#1\hskip-\wd\overbox}$\hskip\wd\overbox}}
\def\overnoarrow#1{\mskip1.5mu\overline{\mskip-1.5mu#1}}
\def\overrightarrow#1{\mskip2mu\vbox{\m@th\ialign{##\crcr
   \rightarrowfill\crcr
   \noalign{\kern-.4pt               
        \kern-\fontdimen22\textfont2 
        \nointerlineskip}
   ${\mskip0mu\hfil\fakeover{#1}\hfil\mskip6mu}$\crcr}}\mskip-2mu}
\def\overleftrightarrow#1{\mskip-3mu\vbox{\m@th\ialign{##\crcr
   \leftarrowfill\hskip-.6em\rightarrowfill\crcr
   \noalign{\kern-.4pt               
        \kern-\fontdimen22\textfont2 
        \nointerlineskip}
   ${\mskip3mu\hfil\fakeover{#1}\hfil\mskip6mu}$\crcr}}\mskip-2mu}

\def\segment#1{{\smash{\overnoarrow{#1}}}}
\makeatother

\begin{document}

\title{The Distortion of a Knotted Curve}

\author{Elizabeth Denne}
\email[Email: ]{edenne@email.smith.edu}
\affiliation{Department of Mathematics \& Statistics, Smith College,
Northampton, MA 01063}

\author{John M. Sullivan}
\email[Email: ]{jms@isama.org}
\affiliation{Institut f\"ur Mathematik, Technische Universit\"at Berlin,
D--10623 Berlin}

\def\overpw#1#2#3{\begin{overpic}[width=#2]{figs/#1}{#3}\end{overpic}}

\begin{abstract}
The distortion of a curve measures the maximum arc/chord length ratio.
Gromov showed any closed curve has distortion at least~$\pi/2$
and asked about the distortion of knots.
Here, we prove that any nontrivial tame knot has distortion at least~$5\pi/3$;
examples show that distortion under~$\upbd$ suffices to build a trefoil knot.
Our argument uses the existence of a shortest essential secant
and a characterization of borderline-essential arcs.
\end{abstract}

\maketitle

Gromov introduced the notion of distortion for curves
as the supremal ratio of arclength to chord length.
(See~\cite{GromovHED}, \cite[p.~114]{GromovFRM} and \cite[pp.~6--9]{GLP}.)
He showed that any closed curve has distortion $\delta\ge\nicefrac\pi2$,
with equality only for a circle.  He then asked whether every knot type
can be built with, say, $\delta\le100$.

\begin{figure*}
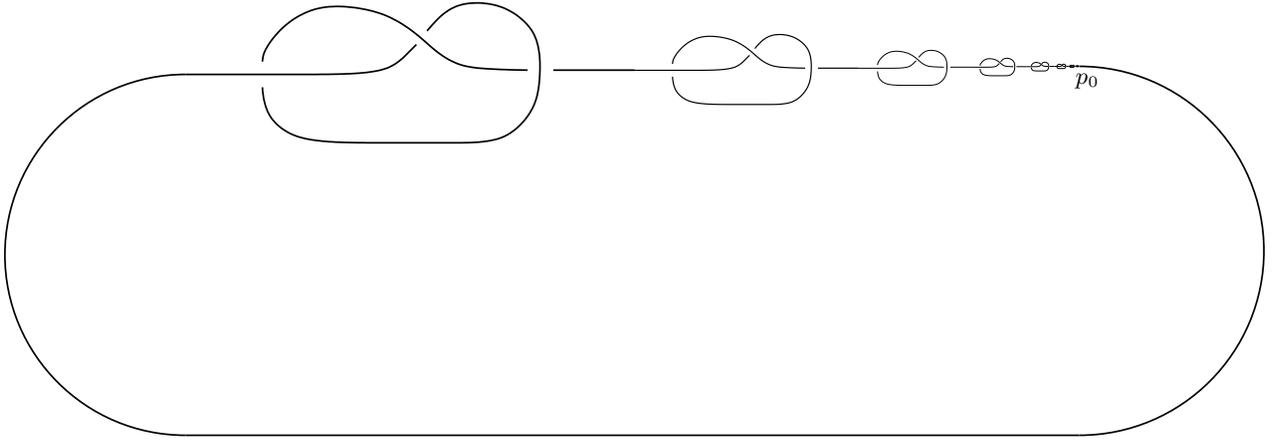

  \overpw{tref-sum}{6.6in}
  { \put(85,28){$p_0$}}
\caption[Wild knot with small distortion]
{A wild knot, the connect sum of infinitely many trefoils,
can be built with distortion less than $10.7$ by repeating
scaled copies of a low-distortion open trefoil.
To ensure that the distortion will be realized within one trefoil,
we merely need to make the copies sufficiently small
compared to the overall loop of the knot
and sufficiently distant from each other.
This knot is smooth except at the one point~$p_0$.}
\label{fig:tref-sum}
\end{figure*}

As Gromov knew, there are infinite families with such a uniform bound.
For instance, an open trefoil (a long knot with straight ends)
can be built with $\delta<10.7$, as follows from
explicit computation for a simple shape.  Then connect sums of
arbitrarily many trefoils---even infinitely many, as in \figr{tref-sum}---can
be built with this same distortion.  (O'Hara~\cite{Ohara1} exhibited
a similar family of prime knots.)

Despite such examples,
many people expect a negative answer to Gromov's question.
We provide a first step in this direction,
namely a lower bound depending on knottedness:
we prove that any nontrivial tame knot has $\delta\ge\newbd$,
more than three times the minimum for an unknot.

To make further progress on the original question, one should try to
bound distortion in terms of some measure of knot complexity.
Examples like \figr{tref-sum} show that crossing number and
even bridge number are too strong: distortion
can stay bounded as they go to infinity.  Perhaps it is worth
investigating hull number \cite{CKKS,Izm-hull}.

Our bound arises from considering \emph{essential secants} of the knot,
a notion introduced by Kuperberg~\cite{Kup}
and developed further in~\cite{DDS}.
There, we used the essential alternating quadrisecants of~\cite{Denne}
to give a good lower bound for the ropelength \cite{GM,CKS2}
of nontrivial knots.

Here, we first show that any knot has a shortest essential secant.
Then we show its endpoints have distortion at least $\newbd\approx5.236$,
using a characterization from~\cite{DDS} of borderline-essential secants.

Our bound is of course not sharp, but numerical simulations~\cite{Mul}
have found a trefoil knot with distortion less than~$\upbd$, so we
are not too far off.  We expect the true minimum distortion for a trefoil
is closer to that upper bound than to our lower bound.
A sharp bound (characterizing that minimum value) would presumably
require a criticality theory for distortion minimizers.
Perhaps this could be developed along the lines of the balance criterion
for (Gehring) ropelength~\cite{CFKSW1},
but the technical difficulties seem formidable.

On the other hand, it is easy to see how to slightly improve our bound
$\delta\ge\newbd$.
Indeed, the circular arc shown in \figr{circ56}
must actually spiral out in the middle (to avoid greater distortion
between~$c$ and~$\0$).
In the first version~\cite{DS-disto-v1} of this paper,
our bounds considered a shortest essential arc;
we used logarithmic spirals to improve an initial
bound $\delta\ge\pi$ to $\delta>3.99$. 
Bereznyak and Svetlov~\cite{BerSve} then obtained $\delta>4.76$
by focusing on a shortest borderline-essential secant and 
making further use of spirals.
We expect that such spirals could improve our bound here
by only a few percent, at the cost of tripling the length of this paper;
thus we have not pursued this idea.

There are easy upper bounds for distortion in terms of
other geometric quantities for space curves.
For instance, an arc of total curvature $\alpha<\pi$
has distortion at most $\sec\nicefrac\alpha2$. (See \cite[\S7]{Sul-FTC}.)
Similarly, a closed curve of ropelength~$R$
has distortion at most $\nicefrac R2$. (This was \cite[Thm.~3]{LSDR}
and also follows easily from \cite[Lem.~3.1]{DDS}.)
But there are no useful bounds the other way: the example
in \figr{tref-sum} has bounded distortion but infinite
total curvature and ropelength, while a steep logarithmic spiral
shows that arcs of infinite total curvature can have distortion
arbitrarily close to~$1$.

This means that a lower bound like ours for the distortion
of a nontrivial knot cannot be based on the known lower bounds
for total curvature~\cite{milnor} or ropelength~\cite{DDS}.
Indeed, before our work here, it remained conceivable that
the infimal distortion of knotted curves was $\nicefrac\pi2$.

An alternative approach might be to consider explicitly
the geometry of curves of small distortion.
A closed plane curve with distortion close to $\nicefrac\pi2$ must be
pointwise close to a round circle~\cite{DEGKR};
we have modified that argument to apply to space curves~\cite{DS-disto-v1}.
But being close to a circle does not preclude being knotted.
We note that the proof looks only at distortion between
opposite points on the curve, and, indeed, any knot type can be realized
so that this restricted distortion is arbitrarily close to~$\nicefrac\pi2$.

\section{Definitions and background}
We deal with oriented, compact, connected curves embedded in $\R^3$.
Such a curve is either an \emph{arc} homeomorphic to an interval,
or a \emph{knot} (a simple closed curve) homeomorphic to a circle.

Two points $p$, $q$ along a knot~$K$ separate~$K$ into two complementary
arcs, $\arc{pq}$ (from~$p$ to~$q$) and $\arc{qp}$.
We let $\len{pq}$ denote the length of~$\arc{pq}$.
Distortion contrasts the shorter arclength distance
$d(p,q):= \min(\len{pq},\len{qp})\le\nicefrac{\lenc{K}}{2}$
with the straight-line (chord) distance $|p-q|$ in~$\R^3$.
(For an arc~$\gamma$, if $p,q$ lie in order along~$\gamma$, then
$d(q,p)=d(p,q):=\len{pq}$ is the length of the subarc~$\arc{pq}$.)

\begin{definition}
The \emph{distortion} between distinct points~$p$ and~$q$
on a curve~$\gamma$ is
$\delta(p,q) := \nicefrac{d(p,q)}{|p-q|}\ge1$.
The \emph{distortion} of~$\gamma$ is the supremum
$\delta(\gamma) := \sup \delta(p,q),$
taken over all pairs of distinct points.
\end{definition}

\begin{figure}
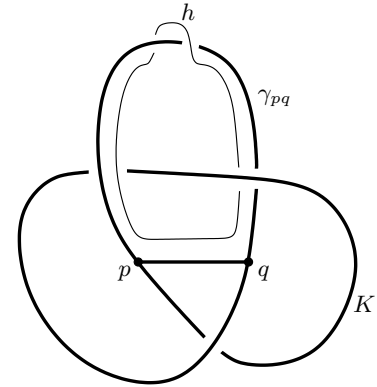

  \overpw{essdef}{1.8in}
  { \put(28,30){$p$} \put(66,30){$q$} \put(92,20){$K$}
    \put(45,100){$h$} \put(66,78){$\arc{pq}$} }
\caption[Essential secants]
{The arc $\arc{pq}$ is essential in the knot~$K$:
the parallel $h$, whose linking number with~$K$ is zero,
is homotopically nontrivial in the knot complement.
This shows there is no disk spanning $\arc{pq}\cup\segment{pq}$
while avoiding~$K$.
In this example, $\arc{qp}$ is also essential, so $\segment{pq}$ is essential.}
\label{fig:essdef}
\end{figure}

Our new bound uses the notion of essential arcs, introduced
in~\cite{DDS} as an extension of ideas of Kuperberg~\cite{Kup}.
Note that generically a knot~$K$ together with a chord $\segment{pq}$
forms a $\theta$-graph in space; being essential is a topological
feature of this knotted graph, as shown in \figr{essdef}.

\begin{definition}
Suppose $\alpha$, $\beta$ and $\gamma$ are interior-disjoint arcs
from~$p$ to~$q$, forming a knotted $\theta$-graph in~$\R^3$.
We say the ordered triple $(\alpha,\beta,\gamma)$ is \emph{essential}
if the loop $\alpha\cup\beta$ bounds no (singular) disk whose interior
is disjoint from the knot $\alpha\cup\gamma$.

Now suppose $K$ is a knot and $p,q\in K$.
If $\segment{pq}$ has no interior intersections with $K$,
we say $\arc{pq}$ is an \emph{essential arc} of~$K$ if
$(\arc{pq},\segment{pq},\arc{qp})$ is essential.
If $\segment{pq}$ does intersect~$K$, we say $\arc{pq}$ is \emph{essential}
if for any $\eps>0$ there is an $\eps$-perturbation~$S$ of~$\segment{pq}$
such that $(\arc{pq},S,\arc{qp})$ is essential.
We say~$\segment{pq}$ is an \emph{essential secant}
if both~$\arc{pq}$ and~$\arc{qp}$ are essential.
\end{definition}

Note that the $\eps$-perturbation ensures that
the set of essential arcs is closed
within the set $(K\times K)\setm\Delta$ of all subarcs.
We say the arc $\arc{pq}$ is \emph{borderline-essential}
if it is in the boundary of the set of essential arcs.  That is,
$\arc{pq}$ is essential, but there are inessential subarcs of~$K$
with endpoints arbitrarily close to~$p$ and~$q$.

The following theorem~\cite[Thm.~7.1]{DDS}
lies at the heart of our distortion bounds.
It describes the special geometric configuration, shown in \figr{bordess},
arising from any borderline-essential arc.
\begin{figure}
  \overpw{bordess}{3.0cm}
  { \put(28,30){$p$} \put(66,89){$q$} \put(44,64){$x$} }
\caption[A borderline-essential arc]
{If the arc $\arc{pq}$ is borderline-essential in the knot~$K$,
\thm{change} gives a point $x\in K\cap\segment{pq}$ for which
$\segment{xp}$ and $\segment{xq}$ are essential.}
\label{fig:bordess}
\end{figure} 

\begin{theorem}\label{thm:change} 
Suppose~$\arc{pq}$ is a borderline-essential subarc of a knot~$K$.
Then the interior of segment~$\segment{pq}$ must intersect~$K$ 
at some point $x\subset\arc{qp}$ for which
the secants~$\segment{xp}$ and~$\segment{xq}$ are both essential.
\qed \end{theorem}

In \cite[Lem.~4.3]{DDS} we showed that the minimum
length of an arc $\arc{ab}\subset\R^n$
staying outside the unit ball~$B_1(\0)$
is $m\big(|a|,|b|,\angle a\0b\big)$,
where for $r,s\ge 1$ and $\theta\in[0,\pi]$ we set
$$m(r,s,\theta) := \begin{cases}
  \sqrt{r^2+s^2-2rs\cos\theta}               &\text{if $\theta\le\theta_0$,} \\
  \sqrt{r^2-1}+\sqrt{s^2-1}+\theta-\theta_0  &\text{if $\theta\ge\theta_0$,}
\end{cases}$$
with $\theta_0=\theta_0(r,s):= \arcsec r + \arcsec s$.
(In the case of plane curves, this can be dated back to~\cite{Kubota}.)

This bound is hard to apply since $m(r,s,\theta)$
is not monotonic in~$r$ and~$s$.
Thus we are led to define $m_1(s,\theta) := \min_{r\ge1} m(r,s,\theta)$,
from which we calculate
$$m_1(s,\theta) = \begin{cases}
  s\sin\theta                     &\text{if $\theta\le \arcsec s$,} \\
  \sqrt{s^2-1}+\theta-\arcsec s   &\text{if $\theta\ge \arcsec s$.}
\end{cases}$$
This function~$m_1$ is continuous, increasing in~$s$ and in~$\theta$,
and concave in~$\theta$.  We have:
\begin{lemma}\label{lem:min}
An arc $\arc{ab}$ staying outside~$B_1(\0)$
has length at least $m_1\big(|b|,\angle a\0b\big)\ge\angle a\0b$.
\qed\end{lemma}

\begin{remark}
For $\theta=\pi$, we are always in the second case
in the definition of~$m_1$, and we have
$m_1(s,\pi)>\sqrt{s^2+1}+\nicefrac\pi2$,
the right-hand side being the length of a curve
that follows a quarter-circle from~$a$ and then
goes straight to~$b$ (cutting into the unit ball).
\end{remark}

\section{Shortest essential arcs and secants}

For a knot of unit thickness, we showed \cite[Lem.~8.1]{DDS}
that essential arcs have length at least~$\pi$,
and essential secants have length at least~$1$.
Here, we show that sufficiently short arcs
and secants of any tame knot are inessential,
and thus that nontrivial tame knots have
shortest essential arcs and secants.

If $K$ is unknotted, any subarc is inessential.
Conversely, Dehn's lemma can be used~\cite[Thm.~5.2]{DDS} to show that
if both~$\arc{pq}$ and~$\arc{qp}$ are inessential (for some $p,q\in K$)
then $K$ is unknotted.  Equivalently, if $K$ is a nontrivial knot
then the complement of any inessential arc is essential.

\begin{lemma}
If $\arc{pq}$ is a borderline-essential subarc of a knot~$K$,
then $\segment{pq}$ is an essential secant.
\end{lemma}

\begin{proof}
Since $\arc{pq}$ is borderline-essential,
there are inessential arcs~$\arc{p'q'}$ converging to~$\arc{pq}$.
Since $K$ (having the essential subarc~$\arc{pq}$) must
be nontrivial, the complements~$\arc{q'p'}$ are essential.
Thus their limit~$\arc{qp}$ is also essential.
\end{proof}

\begin{corollary}\label{cor:esssec}
Given any point~$p$ on a nontrivial knot~$K$, there
is some $q\in K$ such that $\segment{pq}$ is essential.
\end{corollary}

\begin{proof}
Since $K$ is nontrivial, at least some subarcs starting or ending
at~$p$ are essential.  If they all are, then so are all secants from~$p$.
Otherwise there is some borderline-essential arc starting or ending at~$p$.
By the lemma, this gives us an essential secant~$\segment{pq}$.
\end{proof}

\begin{lemma}\label{lem:balliness}
Suppose $K$ is knot and $U$ is a topological ball
such that $K$ intersects~$U$ in a single unknotted arc.
Suppose~$p$ and~$q$ are two points in order along this arc,
and $\beta$ is any arc within~$U$ from~$p$ to~$q$,
disjoint from~$K$. Then $(\arc{pq},\beta,\arc{qp})$ is inessential.
\end{lemma}

\begin{proof}
By definition of an unknotted ball/arc pair, after applying
an ambient homeomorphism we may assume that $U$ is a round
ball and $K\cap U$ a diameter.
Pick any homeomorphism between~$\arc{pq}$ and~$\beta$ (fixing~$p$
and~$q$).  Join all pairs of corresponding points by straight segments;
these fill out a (singular) disk with boundary $\arc{pq}\cup\beta$,
which by convexity stays entirely within~$U$.  The disk avoids~$K$
(except of course for the segment endpoints along~$\arc{pq}$) because
$\beta$ avoids the straight segment $K\cap U$.
\end{proof}

\begin{proposition}\label{prop:shortiness}
Given a knot~$K$ and any locally flat point $a\in K$,
we can find $r>0$ such that
any subarc~$\arc{pq}$ or secant~$\segment{pq}$ of~$K$ 
which lies in the ball $B_r(a)$ is inessential.
\end{proposition}

\begin{proof}
Locally flat means, by definition, that $a$ has a neighborhood~$U$
in which $K\cap U$ is a single unknotted arc.
Choose~$r$ such that $B:=B_r(a)$ is contained in~$U$,
as in \figr{locflat}.
\begin{figure}
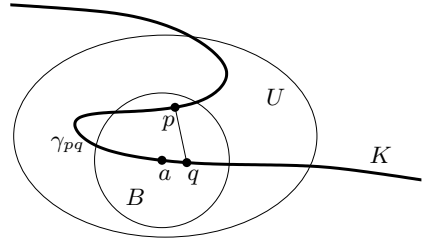

  \overpw{locflat}{55mm}
  { \put(36,14){$a$} \put(43,14){$q$} \put(37,27){$p$}
    \put(87,18){$K$} \put(28,8){$B$} \put(62,32){$U$} \put(10,22){$\arc{pq}$} }
\caption[Short secants are inessential]
{Near a locally flat point $a\in K$, short arcs and secants~$\segment{pq}$
are inessential.}
\label{fig:locflat}
\end{figure}
For any points $p,q\in K\cap B$, the segment~$\segment{pq}$
is contained in~$B$ by convexity, hence in~$U$.
So any sufficiently small perturbation~$S$ of this
segment (as in the definition of essential) stays in~$U$.
Since $K\cap U$ is a single arc, after
switching~$p$ and~$q$ if necessary, we have $\arc{pq}\subset U$.
(If we are proving the first claim, we already know $\arc{pq}\subset B$.)
By \lem{balliness}, $(\arc{pq},S,\arc{qp})$
is inessential, implying by definition that the
subarc~$\arc{pq}$ and the secant~$\segment{pq}$
are inessential.
\end{proof}

\begin{theorem}\label{thm:inessecant}
Given any tame knot~$K$, there exists $\eps>0$ such that
any subarc $\arc{pq}$ of length $\len{pq}<\eps$ is
inessential and any secant~$\segment{pq}$ of length
$|p-q|<\eps$ is inessential.
\end{theorem}

\begin{proof}
Suppose there were sequences $p_n,q_n\in K$ giving
essential arcs or essential secants with length
decreasing to zero.
By compactness of $K\times K$
we can extract a convergent subsequence
$(p_n,q_n)\to (a,a)$.
But the tame knot~$K$ is by definition locally flat at every point $a\in K$.
Choose $r>0$ as in \prop{shortiness} and choose~$n$ large enough that
$\arc{p_nq_n}\subset B_r(a)$.
Then the proposition says $\arc{p_nq_n}$ and $\segment{p_nq_n}$
are inessential, contradicting our choice of $p_n,q_n$.
\end{proof}

\begin{corollary}\label{cor:shortessarc}
Any nontrivial tame knot~$K$ has a shortest essential secant
and a shortest essential subarc.
\end{corollary}

\begin{proof}
Being nontrivial, $K$ does have essential subarcs and secants by \cor{esssec}.
By compactness, a length-minimizing sequence $(p_n,q_n)$ for either case has
a subsequence converging to some $(p,q)\in K\times K$,
and $p\ne q$ by \thm{inessecant}.
Since being essential is a closed condition, this limit arc or secant is
still essential, with minimum length.
\end{proof}


\begin{remark}
A wild knot, even if its distortion is low, can have
arbitrarily short essential arcs, as in the example of \figr{tref-sum}.
For this technical reason, our main theorem will only to tame knots,
even though we expect wild knots must have even greater distortion.
\end{remark}

\section{Distortion bounds}\label{sec:main}

The key to our distortion bounds will be to focus
on a shortest essential secant, as guaranteed by \cor{shortessarc};
we usually rescale so this secant has length~$1$.
Then \thm{change} implies that any
borderline-essential secant has length at least~$2$.
We now bound the length of essential arcs:

\begin{proposition}\label{prop:ess-pi}
Let $K$ be a nontrivial tame knot, scaled so that
a shortest essential secant has length~$1$.
Suppose arc~$\arc{pq}$ is borderline-essential,
and $x\in\segment{pq}\cap K$ is a point as guaranteed by \thm{change}
with~$\segment{xp}$ and~$\segment{xq}$ essential.
Setting $s := \min(|q-x|,2)\in[1,2]$,
we have $\len{pq}\ge m_1(s,\pi)\ge \pi$.
\end{proposition}

\begin{proof}
Translate so that $x$ is the origin~$\0$.
If $\segment{y\0}$ is essential for all $y\in\arc{pq}$,
then by our scaling, $\arc{pq}$ stays outside~$B_1(\0)$.
Thus by \lem{min} and monotonicity of~$m_1$, we get
$\len{pq}\ge m_1(|q|,\pi)\ge m_1(s,\pi)$ as desired.

Otherwise, let $y,z\in\arc{pq}$ be the first and last points
making borderline-essential secants~$\segment{\0y}$ and~$\segment{\0z}$.
By our choice of scaling, $|y|,|z|\ge 2$, and arcs~$\arc{py}$
and~$\arc{zq}$ stay outside $B_1(\0)$.  As in \figr{pcase2},
\begin{figure}
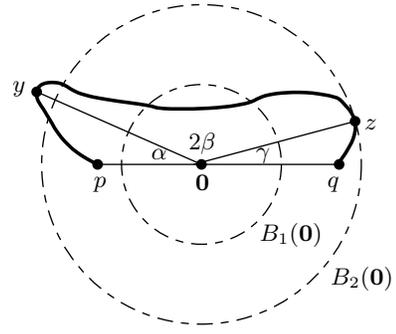

  \overpw{pcase2}{45mm}
  { \put(50,40){$\0$} \put(20,41){$p$} \put(89,41){$q$}
    \put(-4,69){$y$} \put(100,58){$z$} \put(48,52){$2\beta$}
    \put(69,25){$B_1(\0)$} \put(90,12){$B_2(\0)$}
    \put(37,49){$\alpha$} \put(68,49.5){$\gamma$} }
\caption[The case where the arc does not stay essential]
{In the second case in the proof of \prop{ess-pi}, since~$y$ and~$z$
are borderline-essential to~$\0$, they are outside~$B_2(\0)$.
Even though $\arc{yz}$ can go inside~$B_1(\0)$
the total length~$\len{pq}$ in this case is at least $m_1(2,\pi)$.}
\label{fig:pcase2}
\end{figure}
define angles
$$\alpha:=\angle p\0y, \quad
2\beta := \angle y\0z, \quad
\gamma := \angle z\0q,$$
so that $\alpha+2\beta+\gamma=\pi$.  By \lem{min}, we have
$$\len{pq}=\len{py}+\len{yz}+\len{zq}
\ge m_1(2,\alpha)+4\sin\beta+m_1(2,\gamma).$$
By concavity of~$m_1$, for any given~$\alpha+\gamma$, the
sum of the first and last terms is minimized for~$\gamma=0$.
Thus we get $$\len{pq}\ge m_1(2,\pi-2\beta)+4\sin\beta.$$
For $\beta\ge \nicefrac\pi3$, the first case in the
definition of~$m_1$ applies, so
$$\len{pq}\ge 2\sin2\beta+4\sin\beta=4\sin\beta\,(1+\cos\beta)\ge 4.$$
For $\beta\le \nicefrac\pi3$, the second case applies, so
$$\len{pq}\ge \sqrt3+\nicefrac{2\pi}{3} - 2\beta+4\sin\beta\ge
\sqrt3+\nicefrac{2\pi}{3} = m_1(2,\pi).$$
Noting that $4>m_1(2,\pi)\approx 3.826$, we find that
in either case, $\len{pq}\ge m_1(2,\pi)\ge m_1(s,\pi)$.
\end{proof}

\begin{theorem}\label{thm:newbd}
Let $K$ be a nontrivial tame knot, scaled so that
a shortest essential secant has length~$1$.
Suppose $\segment{ab}$ is an essential secant
with length $|a-b|\le 2$.  Then $d(a,b)\ge 2\pi-2\arcsin \nicefrac{|a-b|}2$.
\end{theorem}
\begin{proof}
Switching $a$ and~$b$ if necessary, we may assume $d(a,b)=\len{ab}\le\len{ba}$.
Setting $\phi:=2\arcsin\nicefrac{|a-b|}2\ge|a-b|$, we wish to show that
$\len{ab}\ge2\pi-\phi$.

Let $\arc{ac}\subset\arc{ab}$ be the shortest initial subarc that is essential,
and translate so that the origin $\0\in\segment{ac}\cap K$ is a point as
in \thm{change}.
By \prop{ess-pi}, $\len{ac}\ge m_1(|c|,\pi)\ge \pi$,
so it suffices to show $\len{cb}\ge\pi-\phi$.

For a fixed length $|a-b|$, consider~$\angle a\0b$ as a function
of~$|a|,|b|\ge1$.  It is maximized when $|a|=1=|b|$,
with $\angle a\0b = \phi$.  Thus $\angle c\0b \ge \pi-\phi$.
If $\segment{\0x}$ is essential
for all $x\in \arc{cb}$,
then $\arc{cb}$ remains outside $B_1(\0)$, as in \figr{circ56},
\begin{figure}
  \overpw{circ56}{26mm}
  { \put(41,49){$\0$} \put(96,46){$a$} \put(-2,46){$c$}
    \put(75,8){$b$} 
    \put(52,27){$1$} \put(69,54){$1$} }
\caption[The case where the bound is sharp]
{Since $a$ and~$b$ are outside $B_1(\0)$, we have
$\angle a\0b\le 2\arcsin \nicefrac{|a-b|}2$.
When $\arc{bc}$ stays outside the ball, its length is
at least $\pi-\angle a\0b$.  This figure shows the case $|a-b|=1$
used in \cor{newbd}.}
\label{fig:circ56}
\end{figure}
so $\len{cb}\ge m_1(1,\angle c\0b)=\angle c\0b$
and we are done.

Otherwise, let $x\in\arc{cb}$ be the first point for which $\segment{\0x}$
is borderline-essential, implying that $|x|\ge 2$.
By the triangle inequality
$\len{xb}\ge |x-b| \ge |x-a|-|a-b|$, so
$$\len{cb}\ge \len{cx}+|x-a|-|a-b|.$$
Now set $\theta:=\angle c\0x$ as in \figr{tcase2} and consider two cases.
\begin{figure}
  \overpw{tcase2}{36mm}
  { \put(50,42.5){$\0$} \put(11,47){$c$} \put(84,45){$a$}
    \put(38,43){$\theta$} \put(57,32){$b$} \put(1,11){$x$}
    \put(69,26){$B_1(\0)$} \put(89,12){$B_2(\0)$} }
\caption[The case where the arc does not stay essential]
{If there is a point $x\in\arc{cb}$ with $|x|\ge 2$, then
$\arc{xb}$ can go inside $B_1(\0)$.  We use the bound
$\len{cb} 
\ge |x-c|+\big(|x-a|-|a-b|\big)$.}
\label{fig:tcase2}
\end{figure}

For $\theta\ge\nicefrac\pi2$,
we get $\len{cx}\ge m_1(2,\theta)=\sqrt3+\theta-\nicefrac{\pi}3$,
while $|x-a|\ge 2\sin\theta$ since $|x|\ge 2$.
The concave function $\theta+2\sin\theta$
is minimized at the endpoint $\theta=\pi$, so, as desired,
$$\len{cb}\ge\nicefrac{2\pi}3+\sqrt3-|a-b|>\pi-\phi.$$

For $\theta\le\nicefrac\pi2$, we use $\len{cx}\ge|x-c|$
and consider fixed values of $|c|\ge 1$ and $|x|\ge 2$.
We want to minimize the sum $|x-c|+|x-a|$.
Since $|x-a|$ is increasing in~$|a|$, we may assume $|a|=1$.
Then since $|c|\ge|a|$ and $\theta\le\nicefrac\pi2$, we have
$\angle cx\0 \ge \angle \0xa$.  This means that $|x-c|+|x-a|$
is an increasing function of~$\theta$, minimized at $\theta=0$,
where we have $|x-c|+|x-a|\ge 2-|c|+3$.
Thus $\len{cb}\ge 5-|c|-|a-b|$.
Using the remark after \lem{min},
we have $\len{ac}\ge m_1(|c|,\pi)>|c|+\nicefrac\pi2$.
Thus finally, as desired,
\begin{align*}
\len{ab}
&>\nicefrac\pi2+5-|a-b|>2\pi-\phi. \qedhere
\end{align*}
\end{proof}

\begin{corollary}\label{cor:newbd}
Any nontrivial tame knot has $\delta\ge\newbd$.
\end{corollary}
\begin{proof}
Let $\segment{ab}$ be a shortest essential secant for the knot~$K$,
and scale so that $|a-b|=1$.  Applying the theorem, we get
$\delta(a,b)=d(a,b)\ge\newbd$.
\end{proof}


\bibliography{thick} 
\end{document}